\documentclass [11pt]{article}
\usepackage{amsmath}
\usepackage{amssymb}
\usepackage{amscd}
\usepackage{amsthm}
\usepackage{amsopn}
\usepackage{xspace}
\usepackage{verbatim}
\usepackage{amsmath}
\usepackage{amsfonts}

\newtheorem{theorem}{Theorem}
\newtheorem{lemma}{Lemma}[section]
\newtheorem{proposition}{Proposition}[section]
\newtheorem{definition}{Definition}[section]

\newtheorem{acknowledgment*}{Acknowledgment}
\newtheorem{remark}{Remark}[section]
\numberwithin{equation}{section}

\newcommand{\be}{\begin{equation}}
\newcommand{\ee}{\end{equation}}
\newcommand{\bd}{\begin{displaymath}}
\newcommand{\ed}{\end{displaymath}}

\newcommand{\R}{\mathbb R}

\newcommand{\oz}{\overline{z}}
\renewcommand{\vec}[1]{\boldsymbol{#1}}

\begin{document}
\Large
\begin{center} \bf{Incompressible, quasi-rigid deformations of  2-dimensional  domains }\end{center}
\normalsize
\begin{center} Gershon Wolansky\footnote{ gershonw@math.technion.ac.il}\\
Department of Mathematics \\
Technion, Israel Institute of Technology \\
Haifa 32000, Israel \end{center}
\begin{abstract}
This paper proposes a sensible definition of a deformation metric
between 2-dimensional surfaces obtained from each other by an area
preserving (incompressible) mapping, and  an algorithm for obtaining
this metric, as well as the optimal deformation.
\end{abstract}
\section{Introduction}
 Recently, there is an increasing
interest in the analysis of near-rigid deformation within the
computer vision community, in particular pattern recognition, image
segmentation and face recognition ([B2K], [EK], [MS]).
\par
For example, a variety of objects can be represented as point
clouds. These can be obtained by sampling of the objects in question
by points in some canonical Euclidian space.  One is often presented
with the problem of deciding whether two of these point clouds,
and/or the corresponding underlying objects or manifolds, represent
the same geometric structure (object recognition and
classification).
\par
To quantify the difference between two such clouds  it is natural to
construct smooth domains out of the samplings, and look for a
mapping from one domain  to the other which is "as close to an
isometry as possible". The minimal deviation of these mappings from
a rigid deformation (isometry) can stand as a measure of similarity
between the original objects. \par
 However, this task is very
difficult from a computational point of view, since the set of all
mappings between two domains is very large. On the other hand, it is
sensible to assume that the {\it density} of  the sample points
reflects the true nature of the object. This implies that the {\it
volume elements} associated with these domains, created out of the
samplings, are prescribed. In particular
 we may restrict ourselves to an approximation of an isometry
which preserve  the volume  of  the two domains.
\par
 The
object of this paper is to propose a sensible definition of a
deformation metric between 2-dimensional surfaces obtained from each
other by incompressible (in this case, area preserving)  mapping,
and an algorithm for obtaining this metric, as well as the optimal
deformation.
\par
In section~\ref{example} the problem is formulated for a pair of
flat domains.   Section~\ref{main}  introduces the analytic
conditions for quasi-rigid deformation, for flat domains. In
section~\ref{maingen} the problem is extended to a pair of embedded
surfaces in $\R^3$. In section~\ref{appflow} we use the results of
the previous sections to propose an algorithm, namely a flow which
converges (formally) to a quasi rigid deformation. For the
convenience of the reader we defer all technical proofs to the final
section~\ref{proofs}.
\section{A particular example}\label{example}
 Consider, for example, two
flat domains $\Omega$, $\Omega_1\subset \R^2$, equipped with the
Euclidian metric $e(dx, dy)=dx^2+dy^2$.  If $\Omega_1$ is a rigid
deformation of $\Omega$ then there exists an isometry $\Phi:(\Omega,
e)\rightarrow(\Omega_1,e)$. Locally, it means \be\label{iso}
|\Phi_x|^2\equiv|\Phi_y|^2\equiv 1, \ \ \ \Phi_x\cdot\Phi_y\equiv 0
\ \ \text{on} \ \Omega \ . \ee In general $\Omega$ and $\Omega_1$
are not isometric, and there is no mapping verifying (\ref{iso}).
However,
  $\Phi$ is, {\it by definition}, an isometry between
($\Omega, g_\Phi)$ and $(\Omega_1, e)$ where
\begin{multline}\label{inducedM} g_\Phi(dx, dy):= \left(
                      \begin{array}{cc}
                        dx & dy \\
                      \end{array}
                    \right)|D\Phi|^2\left(
                                          \begin{array}{c}
                                            dx \\
                                            dy \\
                                          \end{array}
                                        \right)   \equiv
 |\Phi_x|^2 dx^2 + 2(\Phi_x\cdot\Phi_y) dx dy + |\Phi_y|^2dy^2 \ .
 \end{multline}
 $$   |D\Phi|^2:= \left(
                             \begin{array}{cc}
                               |\Phi_x|^2 & \Phi_x\cdot\Phi_y \\
                               \Phi_x\cdot\Phi_y & |\Phi_y|^2  \\
                             \end{array}
                           \right)_{(x,y)}  \ . $$
How can we quantify the deviation of $(\Omega_1,e)$   from an
isometric image of $(\Omega,e)$? Certainly, it is related to the
deviation of $g_\Phi$ from the Euclidian metric. Recall that
$g_\Phi$ is represented by a symmetric $2\times 2$ matrix. So, we
consider a real valued function defined on the set of symmetric
$2\times 2$ real  matrices  $S(2;\R)$. Let $h_\alpha: S(2;\R)
\rightarrow \R$, where $\alpha$ is some real parameter (see below),
verifying
$$ h_\alpha(A) = h_\alpha(U A U^*) \ \ \text{for any  } \ U\in O(2;\R) \ \text{and} \ \ A\in S(2;\R) \ , $$
$$ h_\alpha(A)\geq h_\alpha(I)\equiv 0$$
where $I$ is the identity $2\times 2$ matrix, and the equality holds
if and only if $A=I$. Now set \be\label{halpha} H_\alpha(\Omega,
\Omega_1)=\inf_{\Phi\in Diff(\Omega; \Omega_1)} \int_\Omega
h_\alpha(|D\Phi|^2) dxdy \ . \ee By this definition we obtain that
$(\Omega,e)$  and $(\Omega_1,e)$ are isometric if and only if
$H_\alpha(\Omega, \Omega_1)=0$.
\par
What is a natural choice of $h_\alpha$? It must be a function of the
eigenvalues $\lambda_i,\  i=1,2$ of $A$, hence it depends on only 2
arguments, say $tr(A)$ and $det(A)$. Note that \be\label{kdef}
k(A):= tr^2(A)- 4 det(A) = (\lambda_1+\lambda_2)^2
-4\lambda_1\lambda_2 = (\lambda_1-\lambda_2)^2 \geq 0 \ , \ee and
$A=I$ if and only if $k(A)=0$ and $det(A)=1$. This leads us to the
natural
 choice
$$ h_\alpha(A) = k(A)+ \alpha \left( det(A)-1\right)^2 \ , $$
where $\alpha>0$ is a parameter.
\par
In this paper we restrict ourselves to incompressible deformations.
This corresponds to the choice $\alpha=\infty$ which implies the
constraint $det(A)=1$. The adaptation of the quasi-rigid deformation
metric (\ref{halpha}) to the incompressible case is obtained by the
{\it constrained optimization}
 \be\label{hinfty} H_\infty(\Omega, \Omega_1):=\inf_{\Phi\in O(\Omega; \Omega_1)}
\int_\Omega k(|D\Phi|^2) dxdy\ee
 where $O(\Omega; \Omega_1)$ is the set of all {\it area preserving}
 diffeomorphisms  $\Phi:\Omega\rightarrow\Omega_1$, that is,
\be\label{ap}  O(\Omega,\Omega_1):=\left\{  \Phi\in
Diff(\Omega;\Omega_1)  \ ; \ \
|\Phi_x|^2|\Phi_y|^2-(\Phi_x\cdot\Phi_y)^2 \equiv 1 \ \text{on} \
\Omega\right\} \ . \ee

Can we compare the  domains $\Omega, \Omega_1$ using the definition
(\ref{hinfty})? Evidently, $H_\infty(\Omega; \Omega_1)<\infty$ if
and only if the set $O(\Omega; \Omega_1)$ is non empty. By its
definition, a necessary condition is
\be\label{eqa}Area(\Omega)=Area(\Omega_1) \ . \ee By a theorem of
Moser [M] it turns out that condition (\ref{eqa})  is also
sufficient, under rather general conditions.
\par
Assuming (\ref{eqa}), the existence of a minimizer of (\ref{hinfty})
is a much more difficult problem. By the definition (\ref{kdef}) of
$k(\cdot)$ and (\ref{ap}) we pose the following equivalent  open
problem. \vskip.2in\noindent {\bf Open problem}:  {\it  Suppose
(\ref{eqa}).  Is there a minimizer of
  \be\label{hatK}\inf_{\Phi\in O(\Omega;\Omega_1)}  \hat{H}(\Phi) \ \ \ \text{where} \ \ \ \hat{H}(\Phi):= \int_{\Omega}
  tr^2\left(|D\Phi|^2\right)
  dxdy  \  \ ?  \ee }
\section{Main result for flat domain}\label{main}
The set (\ref{ap}) is, formally, an infinite dimensional manifold.
There is an associated right-translation  on this manifold by the
group of area preserving diffeomorphisms \be\label{ap1}
O(\Omega):=\left\{ S \in Diff(\Omega)  \ ; \ \
|S_x|^2|S_y|^2-(S_x\cdot S_y)^2 \equiv 1 \ \text{on} \
\Omega\right\} \ . \ee Indeed, $O(\Omega)$ is a group under
composition, and its action  on $O(\Omega;\Omega_1)$ from the right
is defined by
$$ \Phi\in O(\Omega;\Omega_1) \Longrightarrow \Phi\circ S\in
O(\Omega;\Omega_1) \ . $$ $O(\Omega)$ is a formal Lie group. Its Lie
algebra is given by
   $$ o(\Omega):= \left \{ \vec{v}\in C^\infty
  (\Omega; \R^2) \ \ \nabla\cdot \vec{v} \equiv 0 \ \text{on} \ \Omega \ , \ \vec{v}\cdot
  \hat{n}\equiv 0 \  \text{on} \ \partial\Omega \ . \right\}  \ . $$
  Any smooth flow
$t\rightarrow \Phi^{(t)} \in O(\Omega;\Omega_1)$ is generated by an
orbit $t\rightarrow v^{(t)}\in o(\Omega)$ via \be\label{flow-1}
\frac{d}{dt} \Phi^{(t)} = \vec{v}^{(t)}\left(\Phi^{(t)}\right) \ .
\ee The flow (\ref{flow-1}) can be presented by the Euler equation
 \be\label{flow0} \frac{\partial}{\partial t} \Phi^{(t)} - \vec{v}^{(t)}\cdot
D \Phi^{(t)}=0   \ . \ee Our object is to define $\vec{v}^{(t)}\in
o(\Omega)$
 for which
  $\hat{H}\left( \Phi^{(t)}\right)$ is strictly decreasing, and is
  stationary if and only if $\Phi^{(t)}$ is a critical point of
  $\hat{H}$ in the sense described below.
\par
At this stage it is convenient to consider complex number notation.
Here  we represent $(x,y)\sim x+iy\in
  \mathbb{Z}$. Then
\be\label{compnot} \partial_z=\frac{1}{2}(\partial_x-i\partial_y) \
\ ; \ \
\partial_{\oz}=\frac{1}{2}(\partial_x+i\partial_y) \ , \ dz=dx+idy, \ d\oz= dx-idy , \ dz\wedge d\oz \sim -2i dxdy \ .   \ee
 Lemma~\ref{elementary} below is elementary from the
  divergence theorem and the definition (\ref{compnot}):
 \begin{lemma}\label{elementary}
  If $\Omega$ is simply connected then for any $\vec{v}\in o(\Omega)$ there exists a function $\psi\in
  C^\infty(\Omega)$ so that  $\psi\equiv 0$ on $\partial\Omega$  and
  $\vec{v}= (1/2)(-\psi_y, \psi_x)$ on $\Omega$.
 In particular,  $\vec{v}\sim v_1+iv_2= i\partial_{\oz}\psi$ and
 $$ i\partial_{\oz}\psi d\oz= \frac{\partial\psi}{\partial n} |dz| \
 \text{on} \ \ \partial\Omega$$
 where $\partial\psi/\partial n$ is the (real) outward normal derivative of
 $\psi$ on $\partial\Omega$.
  \end{lemma}
   For  $ \Phi\in
  O(\Omega;\Omega_1)$, let the Hopf function [H]
  \be\label{fdef} f_\Phi(z, \oz):= |\Phi_y|^2 - | \Phi_x|^2
  + 2i \left(\Phi_x\cdot \Phi_y\right) \ .  \ee

  \vskip .3in The main result  is:
  \begin{theorem}\label{mainTh}
  Let $\Phi^0\in O(\Omega;\Omega_1)$. Let $\Phi^{(t)}$ be
  a flow (\ref{flow0}) where $\vec{v}^{(t)}\sim
  i\partial_{\oz}\psi^{(t)}$ satisfying $\psi^{(t)}\in C^\infty(\Omega)$,
  $\psi^{(t)}=0$ on $\partial\Omega$ for any $t\in \R$.
 Then $\Phi^{(t)}\in O(\Omega;\Omega_1)$ for any $t\in\R$ and
\be\label{comway} \frac{d}{dt}\hat{H}\left(\Phi^{(t)}\right)=   4
\mathfrak{Im}\int_\Omega f_{\Phi^{(t)}}\partial^2_{\oz}\psi dxdy \ .
\ee
  \end{theorem}
Now, we are in a position to define quasi-rigid deformation as a
critical point of the functional (\ref{hatK}) on the constraint
manifold (\ref{ap}).
\begin{definition}\label{defcrit}
A mapping $\Phi\in O(\Omega_1;\Omega)$ is quasi-rigid if and only if
$d\hat{H}(\Phi)/dt=0$ for any $\psi\in C^\infty(\Omega)$ verifying
$\psi\equiv 0$ on $\partial\Omega$.
\end{definition}
\begin{proposition}\label{next}
$\Phi:\Omega\rightarrow\Omega_1$ is a quasi rigid deformation  if
and only if (\ref{eqa}) and the following conditions hold:
\begin{description}
\item{i)}\ $\mathfrak{Im}\left(\partial^2_{\oz}f_\Phi\right)=0$ \ on $\Omega$.
\item{ii)} \ $ \mathfrak{Im}\left( f_\Phi \frac{dz}{d\oz}\right) =0$ on $\partial\Omega$. .
\end{description}
\end{proposition}
\begin{remark}
It is interesting to compare the conditions of
Proposition~\ref{next} to the case of {\it harmonic maps} between
two Riemannian  surfaces. In that case the function $f_\Phi$ is
holomorphic, that is, $\partial_{\oz}f_\Phi=0$. This follows from
the stationarity of the Dirichlet functional $\hat{H}$ with respect
to {\em all} parameterizations of the domain $\Omega$. In the case
under consideration, the functional $\hat{H}$ is constrained to the
set $O(\Omega)$, and the criticality condition (i) of the
Proposition is weaker. Note also that (ii) can also be written as
$\mathfrak{Im}(f_\Phi dz^2)=0$ on $\partial\Omega$, while $f_\Phi
dz^2$ is the {\it Hopf differential} (see Definition 1.3.10 of [H]).
\end{remark}

\section{Generalizations}\label{maingen}
Here we generalize the results of section~\ref{main}. Instead of
flat domains $\Omega, \Omega_1$, we consider a pair of smooth,
compact surfaces $\Sigma, \Sigma_1$  embedded in $\R^3$. \par Assume
that $\Sigma, \Sigma_1$ are diffeomorphic   to a canonical domain
$\Delta$. For simplicity we may think about the case  $\Delta:= \{
x^2+y^2\leq 1\}\subset \R^2$. We shall later comment about the case
where $\Omega, \Omega_1$ are embedded manifolds without boundary
(e.g. $\Delta$ is the unit sphere $\mathbb{S}^2$).
\par
Let
$$ {\bf X}:\Delta\rightarrow  \ \  \Sigma= {\bf X}(\Delta)\subset \R^3  \ \ \ , \ \ \ {\bf Y}:\Delta\rightarrow
 \ \ {\bf Y}(\Delta)=\Sigma_1\subset \R^3    \ , $$
 diffeomorphisms from $\Delta$ to
$\Sigma$ and $\Sigma_1$ (respectively), then the mapping
\be\label{map} {\bf Y}\circ {\bf X}^{-1}: \Sigma\rightarrow
\Sigma_1\ee describes a diffeomorphism between the surfaces $\Sigma$
and $\Sigma_1$.
\par
In order to investigate the deviation of $ {\bf Y}\circ {\bf
X}^{-1}$ from an isometry,  we use the mappings ${\bf X}$ and ${\bf
Y}$  to pull back the metrics from $\Sigma$ and $\Sigma_1$
respectively, to $\Delta$.  The geometry of $\Sigma$ is pulled back
to $\Delta$ by the parameterization ${\bf X}$ as
\begin{multline}
g_{\bf X}(dx, dy):= \left(
                      \begin{array}{cc}
                        dx & dy \\
                      \end{array}
                    \right)\left.G_{\bf X}\right|_{(x,y)}\left(
                                          \begin{array}{c}
                                            dx \\
                                            dy \\
                                          \end{array}
                                        \right)   \equiv
 \|{\bf X}_x\|^2 dx^2 + 2\langle{\bf X}_x\cdot{\bf X}_y \rangle dx dy + \|{\bf X}_y\|^2dy^2 \ .
 \end{multline}
  where the standard inner product  $\langle \cdot , \cdot \rangle$,
   $\|\vec{v}\|^2 :=\langle\vec{v}, \vec{v} \rangle$
in $\R^3$ is used, and
 $$\left.G_{\bf X}\right|_{(x,y)}:= \left(
                             \begin{array}{cc}
                               \|{\bf X}_x\|^2 & \langle{\bf X}_x\cdot{\bf X}_y \rangle  \\
                              \langle{\bf X}_x\cdot{\bf X}_y \rangle  & \|{\bf X}_y\|^2  \\
                             \end{array}
                           \right)_{(x,y)} \ . $$
                           The corresponding area element pulled  back from $\Sigma$ to $\Delta$ via ${\bf X}$ is
                                                   $$ \left|{\bf X}_x \wedge {\bf X}_y\right| dxdy =
                                                   det\left( G_{\bf X}\right)dxdy\ . $$
                                                   Similar
                                                   expressions
                                                   hold for the
                                                   metric induced from $\Sigma_1$  on
                                                   $\Delta$ via
                                                   ${\bf Y}$.
\par
Given such  reference parameterizations ${\bf X}$ and ${\bf Y}$, the
set of {\it all} diffeomorphisms from $\Sigma$ to $\Sigma_1$ can be
represented by the set of all diffeomorphisms
$\Phi:\Delta\rightarrow\Delta$ as follows: \be\label{mapPhi} {\bf
Y}\circ\Phi^{-1}\circ  {\bf X}^{-1}: \Sigma\rightarrow \Sigma_1 \ ,
\ \ \ \Phi\in Diff(\Delta) \ .  \ee
  Given  $\Phi\equiv
(\phi_1,\phi_2)\in Diff(\Delta)$ we inquire the deviation of the
mapping (\ref{mapPhi}) from an isometry. The pull-back of the metric
from $\Sigma$ to $\Delta$ by
  ${\bf X}\circ \Phi:= {\bf
X}^\Phi:\Delta\rightarrow\Sigma$  is given by \be\label{isorel}
\left. G_{{\bf X}^\Phi}\right|_{(x,y)}:= D\Phi_{(x,y)} \circ \left.
G_{\bf X}\right|_{\Phi(x,y)}\circ  D^* \Phi_{(x,y)} \  \ee where
$$ D\Phi_{(x,y)} := \left(
              \begin{array}{cc}
               \partial_x \phi_1 & \partial_x\phi_2 \\
                \partial_y\phi_1 & \partial_y\phi_2 \\
              \end{array}
            \right)  \  \ , \ \ D^*\Phi_{(x,y)} := \left(
              \begin{array}{cc}
               \partial_x \phi_1 & \partial_y\phi_1 \\
                \partial_x\phi_2 & \partial_y\phi_2 \\
              \end{array}
            \right)   \ . $$
Likewise, the area element from $\Sigma$ is pulled back to this on
$\Delta$ by
$$ \left|{\bf X}^\Phi_x \wedge {\bf X}^\Phi_y\right|_{(x,y)} dxdy = det\left( \left. G_{\bf
X}\right|_{\Phi(x,y)} \right)det^2\left( D\Phi_{(x,y)}\right)
 dxdy
 \ . $$
  In particular, for any $\Phi\in Diff(\Delta)$, the surface
 area $|\Sigma|$ of $\Sigma$ is given by
  $$ |\Sigma|\equiv \int_\Delta det\left( \left. G_{\bf
X}\right|_{\Phi(x,y)} \right)det^2\left( D\Phi_{(x,y)}\right)
 dxdy \ . $$
 and is independent of $\Phi$. Lemma~\ref{Surlem} below follows from
 definition:
 \begin{lemma}\label{Surlem}
 The diffeomorphism (\ref{mapPhi}) is an isometry between $\Sigma$
 and $\Sigma_1$ if and only if
 \be\label{cond1} D\Phi_{(x,y)} \circ \left.
G_{\bf X}\right|_{\Phi(x,y)}\circ  D^* \Phi_{(x,y)} = \left.G_{\bf
Y}\right|_{(x,y)} \  \ \text{on} \ \ \Delta \ . \ee Likewise,
(\ref{mapPhi}) is area preserving if and only if \be\label{cond2}
det\left( \left. G_{\bf X}\right|_{\Phi(x,y)} \right)det^2\left(
D\Phi_{(x,y)}\right)= det\left( \left. G_{\bf Y}\right|_{(x,y)}
\right) \  \text{on} \ \Delta \ . \ee
 \end{lemma}
We now choose  a particular  reference parameterizations ${\bf X}$,
${\bf Y}$ as follows. Recall that the  {\it Uniformisation Theorem}
[FK] implies that  there exist {\it conformal} parameterizations of
these surfaces. We may assume, therefore, that ${\bf X}$ (res. ${\bf
Y}$) are conformal parameterizations to $\Sigma$, (res. $\Sigma_1$).
This means  \be\label{confpar}G_{\bf X}\equiv \mu \left(
                                                      \begin{array}{cc}
                                                        1 & 0 \\
                                                        0 & 1 \\
                                                      \end{array}
                                                    \right) \ \ \ , \ \ \  G_{\bf Y}\equiv \eta  \left(
                                                      \begin{array}{cc}
                                                        1 & 0 \\
                                                        0 & 1 \\
                                                      \end{array}
                                                    \right) \  \ \text{on} \
                                                    \Delta \ ,
 \ee
where $\mu^2$ (res. $\eta^2$) are the area densities on $\Delta$
associated with the conformal parameterizations of $\Sigma$ (res.
$\Sigma_1$). With these particular parameterizations, condition
(\ref{cond1}) for
 isometry of the mapping (\ref{mapPhi}) is reduced to
 \be\label{cond2*} \mu(\Phi(x,y)) \left| D\Phi\right|^2_{(x,y)} = \eta(x,y) \left(
                                                                 \begin{array}{cc}
                                                                   1 & 0 \\
                                                                   0 & 1 \\
                                                                 \end{array}
                                                               \right)
 \ \
\text{on} \ \ \Delta \ , \ee where $\left| D\Phi\right|^2=
D\Phi\circ D^*\Phi$. Likewise, (\ref{mapPhi}) is area preserving if
and only if
$$ \mu(\Phi(x,y))
det\left( D\Phi_{(x,y)}\right)= \eta(x,y) \ \ \ \text{on} \  \Delta
\ .
$$
Let us consider the set
 \be\label{OSigdef} O(\Sigma, \Sigma_1):= \left\{ \Phi\in Diff(\Delta) \ ; \ \  \mu(\Phi(x,y))
det\left( D\Phi_{(x,y)}\right)= \eta(x,y) \ \ \ \forall \ (x,y)\in
\Delta    \right\} \ . \ee  Again, the result of Moser [M] implies
that $O(\Sigma, \Sigma_1)$ is not empty, provided
\be\label{eqar}|\Sigma|=|\Sigma_1| \ . \ee

Now, we recall (\ref{kdef}). It implies \be\label{keq1} k\left(
D\Phi_{(x,y)}\left.G_{\bf X}\right|_{\Phi(x,y)}
D^*\Phi_{(x,y)}\right) \geq 0\ee  and, by Lemma~\ref{Surlem}, it
follows that $\Sigma$, $\Sigma_1$ are isometric if and only if
\be\label{keq2}k\left( D\Phi_{(x,y)}\left.G_{\bf
X}\right|_{\Phi(x,y)} D^*\Phi_{(x,y)}\right) \equiv 0 \ \ \text{on}
\  \Delta \ . \ee Under the assumption (\ref{confpar}) we obtain
$$  k\left( D\Phi_{(x,y)}\left.G_{\bf X}\right|_{\Phi(x,y)}
D^*\Phi_{(x,y)}\right) = \mu^2(\Phi(x,y))tr\left(|D\Phi|^2_{(x,y)}
\right) - 4\mu^2(\Phi(x,y)) det\left( |D\Phi|_{(x,y)}^2\right) \ \ .
$$
For  $\Phi\in O(\Sigma, \Sigma_1)$ it follows \be\label{keq3}k\left(
D\Phi_{(x,y)}\left.G_{\bf X}\right|_{\Phi(x,y)}
D^*\Phi_{(x,y)}\right) = \mu^2(\Phi(x,y))tr\left(|D\Phi|^2_{(x,y)}
\right) - 4\eta^2(x,y)  \ . \ee Using (\ref{keq1},\ref{keq2}) and
(\ref{keq3}) we obtain
\begin{lemma} Assume (\ref{eqar}). If ${\bf X}, {\bf Y}$ are conformal parameterizations and
 $\Phi\in O(\Sigma, \Sigma_1)$  then  \be\label{lemdef} \hat{H}_\Delta(\Phi):=
\int_\Delta \mu^2(\Phi) tr\left(|D\Phi|^2\right) dxdy \geq 4
|\Sigma_1| \ . \ee  The equality in (\ref{lemdef}) holds if and only
if $\Sigma$ and $\Sigma_1$ are isometric.
\end{lemma}
We are now at a position to generalize (\ref{hatK}): Calculate
\be\label{hatKnew}\inf_{\Phi\in O(\Sigma, \Sigma_1)}
\hat{H}_\Delta(\Phi)   \ee and find the optimal mapping $\Phi$ (if
exists).

  As in section~\ref{main}, the right action on the "manifold" $O(\Sigma, \Sigma_1)$
  is given by the Lie group of $\eta$ preserving diffeomorphisms on $\Delta$:
\be\label{OX} O_\eta(\Delta):= \left\{ S\in Diff(\Delta) \ ; \ \
\eta(S(x,y)) det\left( DS_{(x,y)}\right)= \eta(x,y) \ \ \ \forall \
(x,y)\in \Delta \ .  \ \right\} \ . \ee
\begin{lemma} \label{oOSigma}The set $O_\eta(\Delta)$ is a group under compositions.
For any $\Phi\in O(\Sigma, \Sigma_1)$ and $S\in O_\eta(\Delta)$,
$\Phi\circ S\in O(\Sigma, \Sigma_1)$. In particular,
$O_\eta(\Delta)$ acts on $O(\Sigma, \Sigma_1)$ by
right-translations.
\end{lemma}

 The Lie algebra associated with $O_\eta(\Delta)$ is
  \be\label{oXdef} o_\eta:= \left \{ \vec{v}\in C^\infty
  (\Delta; \R^2) \ \ \nabla\cdot\left( \eta\vec{v}\right) \equiv 0 \ \text{on} \ \Delta \ , \ \vec{v}\cdot
  \hat{n}\equiv 0 \  \text{on} \ \partial\Delta \ . \right\}  \ .
  \ee
  In complex notation (\ref{compnot}), Lemma~\ref{elementary}
  implies that
  $$ \vec{v}\in o_\eta \ \text{if and only if} \ \  \vec{v}\sim v_1+iv_2 =
 i \eta^{-1}\partial_{\oz}\psi $$
  where $\psi\in C^\infty(\Delta)$ and $\psi\equiv 0$ on
  $\partial\Delta$.
  \par
As in the case of planar domains discussed in section~\ref{main}, we
consider now the flow $\Phi^{(t)} \in O(\Sigma, \Sigma_1)$ generated
by the Euler equation
 \be\label{flownew} \frac{\partial}{\partial t} \Phi^{(t)} - \vec{v}^{(t)}\cdot
D \Phi^{(t)}=0   \  \ , \ \ \vec{v}^{(t)}\in o_\eta \ ,  \ee and the
Hopf function $f_\Phi$ defined, as in (\ref{fdef}), on $\Delta$. The
generalization of Theorem~\ref{mainTh} is now formulated as
\begin{theorem}\label{mainThSur}
  Let $\Phi^0\in O(\Sigma, \Sigma_1)$. Let $\Phi^{(t)}$ be
  a flow (\ref{flownew}) where $\vec{v}^{(t)}\sim
  i\eta^{-1}\partial_{\oz}\psi^{(t)}$ satisfying $\psi^{(t)}\in C^\infty(\Omega)$,
  $\psi^{(t)}=0$ on $\partial\Omega$ for any $t\in \R$.
 Then $\Phi^{(t)}\in O(\Sigma, \Sigma_1)$ for any $t\in\R$ and
 \be \label{finalR}\frac{d}{dt}\hat{H}_\Delta\left(\Phi^{(t)}\right)
= 4\mathfrak{Im}\int_\Delta
\mu^2\left(\Phi^{(t)}\right)f_{\Phi}\partial_{\oz}\left( \eta^{-1}
\partial_{\oz}\psi\right)dxdy  \ .
\ee
  \end{theorem}
  Now, we generalize Definition~\ref{defcrit} and
  Proposition~\ref{next} as follows:
  \begin{definition}\label{defcritDelta}
A mapping ${\bf Q}:\Sigma\rightarrow\Sigma_1$ is quasi-rigid if and
only if it can be decomposed as \be\label{defxyphi} {\bf Q}=  \bf
{Y}\Phi^{-1}{\bf X}^{-1}\ee where ${\bf X}:\Delta\rightarrow
\Sigma$, ${\bf Y}:\Delta\rightarrow \Sigma_1$ are conformal
diffeomorphisms and $\Phi\in O_\eta(\Delta)$ verifying
$d\hat{H}_\Delta (\Phi)/dt=0$ for any $\psi\in C^\infty(\Delta)$ for
which $\psi\equiv 0$ on $\partial\Delta$.
\end{definition}
\begin{proposition}\label{nextDelta} ${\bf Q}$ given by
(\ref{defxyphi}) is a quasi rigid deformation  verifying
(\ref{confpar}) if and only if the corresponding $\Phi\in
O_\eta(\Delta)$ verifies
\begin{description}
\item{i)}\ $\mathfrak{Im}\partial_{\oz}\left\{ \partial_{\oz}\left( \mu^2(\Phi) f_{\Phi} \right)\eta^{-1}\right\}=0$  \ on $\Delta$.
\item{ii)} \ $ \mathfrak{Im}\left( f_\Phi \frac{dz}{d\oz}\right) =0$ on $\partial\Delta$.
\end{description}
\end{proposition}
  \section{Applications}\label{appflow}
  Theorem~\ref{mainTh} in section~\ref{main} suggests an algorithm
  for calculating a quasi rigid deformations from a flat domain
  $\Omega\subset\R^2$ to another $\Omega_1\subset\R^2$ where (\ref{eqa}) is assumed.
  \begin{theorem}\label{flowTh1}
  Let $\Phi^{(0)}\in O(\Omega)$. Define the flow $\Phi^{(t)}$ in
  $O(\Omega)$ given by the Euler equation (\ref{flow0}) for $t\geq 0$ and
  $\vec{v}^{(t)}\sim i\partial_{\oz}\psi^{(t)}\in o(\Omega)$ is defined as follows:  $\psi^{(t)}=\psi^{(t)}_0+
\psi^{(t)}_1+\psi^{(t)}_2$ where \be\label{maineq} \psi^{(t)}_0:=
-\mathfrak{Im}\left(\partial^2_{\oz}f_{\Phi^{(t)}}\right) \ee on
$\Omega$ and
$$ \int_\Omega \psi^{(t)}_1\psi^{(t)}_0 dxdy=  \int_\Omega
\psi^{(t)}_2\psi^{(t)}_0=0$$ as well as
$$  \psi^{(t)}_1= -\psi^{(t)}_0 \ \text{and} \ \  \partial_z\psi^{(t)}_1=-\partial_z
\psi^{(t)}_0 \ \ \text{on} \ \partial\Omega $$ and $$ \psi^{(t)}_2=0
\ \ \ , \ \ \frac{\partial\psi^{(t)}_2}{\partial n} = -
\mathfrak{Im}\left( f_\Phi \frac{dz}{d\oz}\right) \ \ \text{on} \ \
\partial\Omega \ , $$
then the flow (\ref{flow0}) satisfies $$
\frac{d\hat{H}(\Phi^{(t)})}{dt} =-
2\int_\Omega\int\left|\mathfrak{Im}\left( \partial^2_{\oz}
f_{\Phi^{(t)}} \right)\right|^2 dxdy-
\int_{\partial\Omega}\left|\mathfrak{Im}\left(
f_{\Phi^{(t)}}\frac{dz}{d\oz}\right)\right|^2
 |dz| \ .
$$
In particular, $d\hat{H}(\Phi)/dt\leq 0$ and $\Phi$ is a steady
state of (\ref{flow0}) if and only if it is a quasi-rigid
deformation.
\end{theorem}
Similarly,
  Theorem~\ref{mainThSur} in section~\ref{maingen} suggests an algorithm
  for calculating a quasi rigid deformations from an embedded
  surface $\Sigma\subset \R^3$ to another $\Sigma_1\subset \R^3$:
 \begin{theorem}\label{flowTh2}
 Let $\Sigma\subset \R^3, \Sigma_1\subset \R^3$ and a pair of
 conformal parameterizations ${\bf X}, {\bf Y}$ verifying
 (\ref{confpar}). Let $\Phi^{(0)}\in O_\eta(\Delta)$. Define the flow $\Phi^{(t)}$ in
  $O_\eta(\Delta)$ given by the Euler equation (\ref{flownew}) for $t\geq 0$ and
  $\vec{v}^{(t)}\sim i\partial_{\oz}\psi^{(t)}\in o_\eta$ is defined
  as in Theorem~\ref{appmainTh}, where (\ref{maineq}) is replaced by
\be\label{maineqalt} \psi_0^{(t)}=
-\mathfrak{Im}\partial_{\oz}\left\{
\partial_{\oz}\left( \mu^2(\Phi^{(t)}) f_{\Phi^{(t)}}
\right)\eta^{-1}\right\} \ . \ee then the flow (\ref{flownew})
satisfies

$$ \frac{d\hat{H}_\Delta(\Phi^{(t)})}{dt} =-2\int_\Delta\left|\mathfrak{Im} \partial_{\oz}\left\{
\partial_{\oz}\left( \mu^2(\Phi^{(t)}) f_{\Phi^{(t)}}
\right)\eta^{-1}\right\}\right|^2 dxdy -
\int_{\partial\Delta}\eta^{-1}
\mu^2(\Phi^{(t)})\left|\mathfrak{Im}\left(
f_{\Phi^{(t)}}\frac{dz}{d\oz}\right)\right|^2
 |dz| \ .
$$
In particular, $d\hat{H}(\Phi)/dt\leq 0$ and $\Phi$ is a steady
state of (\ref{flownew}) if and only if ${\bf Y} \Phi^{-1}{\bf
X}^{-1}$ is a quasi-rigid deformation of $\Sigma$ into $\Sigma_1$.
\end{theorem}
\begin{remark}
In case the surfaces $\Sigma, \Sigma_1$ are closed (e.g. both
homeomorphic to the sphere $ \mathbb{S}^2$), then
$\psi^{(t)}=\psi_0^{(t)}$ as given in (\ref{maineqalt}), and there
is no need of the component $\psi^{(t)}_1, \psi^{(t)}_2$ which take
care of the boundary.
\end{remark}
\section{Proofs}\label{proofs}
  \begin{proof}(of Theorem~\ref{mainTh})
 Consider an orbit $\Phi^{(t)}$ induced by $S^{(t)}\in O(\Omega):= O(\Omega; \Omega)$ via
  $$ \Phi^{(t)} = \Phi^{(0)}\circ S^{(t)} \ . $$
A tangent of this orbit is given by the left representation
 \be\label{flow}  \dot{S} = \vec{v}\circ S  \ , \vec{v}\in
o(\Omega) \ , \ \ S^{(0)}=\bf{I} \ .
 \ee
 Since $\nabla\cdot{\vec v}\equiv 0$ by definition it follows by  (\ref{inducedM})that
 $det (DS^{(t)})= det\left(G_{S^{(t)}}\right)\equiv 1$ along this orbit. In particular.
$$ det\left(G_{\Phi^{(0)}\circ S^{(t)}_{(z)}}\right)
= det\left(G_{\Phi^{(0)}}\right)_{{(S^{(t)})}}\cdot
det\left(G_{S^{(t)}}\right)_{(z)} = 1$$ so
 $\Phi^{(t)}\equiv
 \Phi^{(0)}\circ S^{(t)}\in O(\Omega;\Omega_1)$.  In addition
   \be\label{Dflow}  D\dot{S} = \left[D\vec{v}\circ S\right] DS  \ .
 \ee
 In the left representation $\Phi^{(t+\tau)}=\Phi^{(t)}\circ
 S^{(\tau)}$, $S^{(\tau)}(x)= x + \tau \vec{v}(x) + O(\tau^2)$, so
\be\label{expand} D\Phi^{(t+\tau)} = \left[ D\Phi^{(t)}\circ
S^{(\tau)}\right] D S^{(\tau)}= D\Phi^{(t)}+ \tau D\Phi^{(t)}
D\vec{v} + \tau D^2\Phi^{(t)} \vec{v}+O(\tau^2) \ , \ee We obtain
from (\ref{hatK}) and (\ref{expand})
\begin{multline}\label{expand2}\hat{H}\left(\Phi^{(t+\tau)}\right)=\hat{H}\left(\Phi^{(t)}\right)+\tau
\int_\Omega tr \left( D(\Phi^{(t)})\left( D\vec{v}+
D^*\vec{v}\right) D^*\Phi^{(t)}\right) dxdy \\ + \tau \int_\Omega tr
\left( \left[D^2(\Phi^{(t)})D^*\Phi^{(t)}+ D\Phi^{(t)} D
D^*\Phi^{(t)}\right]\vec{v}\right) dxdy + O(\tau^2) \ .
\end{multline}
However
$$ \left[D^2(\Phi^{(t)})D^*\Phi^{(t)}+ D\Phi^{(t)} D
D^*\Phi^{(t)}\right]= D \left|D\Phi^{(t)}\right|^2$$ and integration
by parts implies \be\label{expand3} \int_\Omega D \left\{
tr\left(\left|D\Phi^{(t)}\right|^2\right)\right\}\vec{v}dxdy=-\int_\Omega
tr\left(\left|D\Phi^{(t)}\right|^2\right)(\nabla\cdot\vec{v})dxdy=0
\ , \ee since $\vec{v}$ is divergence free. From (\ref{expand2},
\ref{expand3}) we obtain \be\label{mainL}
 \frac{d}{dt}\hat{H}\left(\Phi^{(t)}\right)=  \int_\Omega
tr\left( D(\Phi^{(t)})\left[ D\vec{v} + D^*\vec{v}\right]
D^*(\Phi^{(t)})\right) dxdy   \ee
 Next, note that if $\Omega$ is
simply connected, then any smooth divergence free vector field can
be written as \be \label{psidef}\vec{v}= \nabla^\bot\psi:= (-\psi_y,
\psi_x), \ \ \psi\in C^\infty(\Omega) \ \ , \psi\equiv 0 \ \
\text{on} \
\partial\Omega \ .
\ee It follows that \be\label{DvD} \frac{1}{2}\left(D\vec{v}+
D^*\vec{v}\right) =\frac{1}{2} {\bf S}(\psi) + {\bf U}(\psi)\ee
                            where, using $\Box:=
                            \partial^2_x-\partial^2_y$,
\be\label{defsu} {\bf S}(\psi):= \left(
                     \begin{array}{cc}
                       0 & \Box\psi \\
                       \Box\psi & 0 \\
                     \end{array}
                   \right) \ \ \ , \ \ \  {\bf U}(\psi):= \left(
                     \begin{array}{cc}
                       -\psi_{x,y} & 0 \\
                       0 & \psi_{x,y} \\
                     \end{array}
                   \right) \ . \ee
                   Note that, with $\Phi^{(t)}:= (\phi_1^{(t)},
                   \phi_2^{(t)})$,
                   $$ D\Phi^{(t)} D^*\Phi^{(t)}\equiv \left|
                   D\Phi^{(t)}\right|^2 = \left(
                                            \begin{array}{cc}
                                              \left|\Phi_x^{(t)}\right|^2 &  \Phi_x^{(t)} \cdot\Phi_y^{(t)} \\
                                               \Phi_x^{(t)} \cdot\Phi_y^{(t)} &   \left|\Phi_y^{(t)}\right|^2 \\
                                            \end{array}
                                          \right) \ , $$
                                          so
\be\label{sqsi} tr\left( D\Phi^{(t)} {\bf S}(\psi)
D^*\Phi^{(t)}\right)= 2
\left(\Phi_x^{(t)}\cdot\Phi_y^{(t)}\right)\Box\psi \ee while
\be\label{crosspsi} tr\left( D\Phi^{(t)} {\bf U}(\psi)
D^*\Phi^{(t)}\right)=
\left(\left|\Phi_y^{(t)}\right|^2-\left|\Phi_x^{(t)}\right|^2\right)\psi_{xy}
\ . \ee Hence, \be\label{treq} tr\left( D(\Phi^{(t)})\left[ D\vec{v}
+ D^*\vec{v}\right] D^*(\Phi^{(t)})\right) =
\left(\Phi_x^{(t)}\cdot\Phi_y^{(t)}\right)\Box\psi +
\left(\left|\Phi_y^{(t)}\right|^2-\left|\Phi_x^{(t)}\right|^2\right)\psi_{xy}
\ . \ee We now return to the complex notation. From (\ref{compnot},
\ref{fdef}), \be\label{treq1}
\left(\Phi_x^{(t)}\cdot\Phi_y^{(t)}\right)\Box\psi +
\left(\left|\Phi_y^{(t)}\right|^2-\left|\Phi_x^{(t)}\right|^2\right)\psi_{xy}=
4 \mathfrak{Im}\left(  f_{\Phi^{(t)}} \partial^2_{\oz}\psi  \right)
\ , \ee hence, by (\ref{mainL}) and (\ref{treq}, \ref{treq1}), $$
\frac{d}{dt}\hat{H}\left(\Phi^{(t)}\right)=  - 2\mathfrak{Re}\left(
\int_\Omega f_{\Phi^{(t)}}\partial^2_{\oz}\psi dz \wedge
d\oz\right)\equiv 4 \mathfrak{Im}\left( \int_\Omega
f_{\Phi^{(t)}}\partial^2_{\oz}\psi dxdy\right) \ .$$
\end{proof}
\begin{proof} (of Proposition~\ref{next}): \\
We imply integration by parts on (\ref{comway}) to obtain
$$ \frac{d}{dt}\hat{H}\left(\Phi^{(t)}\right) = 4\mathfrak{Im}\left(\int_\Omega\int \psi\partial^2_{\oz}
f_{\Phi^{(t)}}  dxdy \right)+
2\mathfrak{Re}\left(\int_{\partial\Omega}\left(\psi
\partial_{\oz}f_{\Phi^{(t)}}  -
  f_{\Phi^{(t)}}\partial_{\oz}\psi \right)dz\right) \ .
$$
However, $\psi\equiv 0$ on $\partial\Omega$ by (\ref{psidef}) so
\be\label{finalRR} \frac{d}{dt}\hat{H}\left(\Phi^{(t)}\right) =
4\mathfrak{Im}\int_\Omega\int  \psi\partial^2_{\oz} f_{\Phi^{(t)}}
dxdy -
2\mathfrak{Re}\int_{\partial\Omega}f_{\Phi^{(t)}}\partial_{\oz}\psi
 dz \ .
\ee Part (i) follows from the first integral (\ref{finalRR}), since
$\psi$ is arbitrary in the interioir of $\Omega$. Since $\psi\equiv
0$ on $\partial\Omega$, it follows by Corollary~\ref{elementary}
that $i\partial_{\oz}\psi d\oz$ is real valued on $\partial\Omega$.
Hence \be\label{finalbc} \mathfrak{Re}\left(
f_\Phi\partial_{\oz}\psi dz\right)= - \mathfrak{Im}\left(
f_\Phi\frac{dz}{d\oz}i\partial_{\oz}\psi
d\oz\right)=\pm\partial_n\psi \mathfrak{Im}\left(
f_\Phi\frac{dz}{d\oz}\right) \ . \ee Since $\partial_n\psi$ is
arbitrary on $\partial\Omega$ we obtain part (ii).
\end{proof}
\begin{proof} (of Lemma~\ref{oOSigma}): \
Let $S, S_1\in O_\eta(\Delta)$. Then $D(S\circ S_1) =
DS_{(S_1)}\circ D S_1$. Hence $ det\left(D(S\circ S_1)\right) =
det\left( DS\right)_{(S_1)} det\left( DS_1\right)$. By definition,
$$ det(DS)= \eta/ \eta(S) \ \ ; \ \  det\left(DS\right)_{(S_1)}= \eta(S_1)/ \eta(S\circ
S_1)$$ so  $det\left(D(S\circ S_1)\right) =\eta/ \eta(S\circ S_1)$,
so $S\circ S_1\in O_\eta(\Delta)$. The same argument holds for
$\Phi\circ S$ where $\Phi\in O(\Sigma, \Sigma_1)$.
\end{proof}
\begin{proof} (of Theorem~\ref{mainThSur}):
\\
We repeat the proof of Theorem~\ref{main} up to (\ref{expand2}).
Here (\ref{expand2}). is replaced by
\begin{multline}\label{expand2*}\hat{H}_\Delta\left(\Phi^{(t+\tau)}\right)=\hat{H}_\Delta\left(\Phi^{(t)}\right)+\tau
\int_\Delta \mu^2\left(\Phi^{(t)}\right)tr \left(
D(\Phi^{(t)})\left( D\vec{v}+ D^*\vec{v}\right) D^*\Phi^{(t)}\right)
dxdy \\ + \tau \int_\Delta \mu^2\left(\Phi^{(t)}\right) tr \left(
\left[D^2(\Phi^{(t)})D^*\Phi^{(t)}+ D\Phi^{(t)} D
D^*\Phi^{(t)}\right]\vec{v}\right) dxdy +  \\ \tau  \int_\Delta
\left(\nabla\mu^2_{(\Phi^{(t)})}\cdot \vec{v}\right)tr \left(|
D\Phi^{(t)}|^2\right) dxdy + O(\tau^2) \ .
\end{multline}
As in (\ref{expand3}), the third term on the right of
(\ref{expand2*}) is reduced to \be\label{expand3**} \int_\Omega
\mu^2(\Phi^{(t)})D \left\{
tr\left(\left|D\Phi^{(t)}\right|^2\right)\right\}\vec{v}dxdy=-\int_\Omega
tr\left(\left|D\Phi^{(t)}\right|^2\right)\nabla\cdot(\mu^2(\Phi^{(t)})\vec{v})dxdy
\ . \ee Together with the forth term in (\ref{expand2*}) we obtain
 \begin{multline}\label{expand2**}\hat{H}_\Delta\left(\Phi^{(t+\tau)}\right)=\hat{H}_\Delta\left(\Phi^{(t)}\right)+\tau
\int_\Delta \mu^2\left(\Phi^{(t)}\right)tr \left(
D(\Phi^{(t)})\left( D\vec{v}+ D^*\vec{v}\right) D^*\Phi^{(t)}\right)
dxdy \\ - \tau \int_\Delta \mu^2_{(\Phi^{(t)})}
(\nabla\cdot\vec{v})tr \left(| D\Phi^{(t)}|^2\right) dxdy +
O(\tau^2) \ .
\end{multline}
Since $\vec{v}\in o_\eta$ it follows from (\ref{oXdef}) and
(\ref{expand2**})
\begin{multline}\hat{H}_\Delta\left(\Phi^{(t+\tau)}\right)=\hat{H}_\Delta\left(\Phi^{(t)}\right)+\tau
\int_\Delta \mu^2\left(\Phi^{(t)}\right)tr \left(
D(\Phi^{(t)})\left( D\vec{v}+ D^*\vec{v}\right) D^*\Phi^{(t)}\right)
dxdy \\ + \tau \int_\Delta \mu^2(\Phi^{(t)})\eta^{-1}( \vec{v}\cdot
\nabla\eta)tr \left(| D\Phi^{(t)}|^2\right) dxdy + O(\tau^2) \ .
\end{multline}
We now use $\vec{v}\sim i\eta^{-1}\partial_{\oz}\psi$ and
(\ref{psidef}-\ref{DvD}) to obtain  $$ D\vec{v}+ D^*\vec{v} =
\eta^{-1}{\bf S}(\psi) +2\eta^{-1}{\bf U}(\psi)+ 2\eta^{-2}\left(
                                  \begin{array}{cc}
                                    -\eta_x\psi_y &  \frac{\eta_x\psi_x-\eta_y\psi_y}{2} \\
                                     \frac{\eta_x\psi_x-\eta_y\psi_y}{2} & \eta_y\psi_x \\
                                  \end{array}
                                \right) \ .
 $$
By direct calculation  \begin{multline}2\mu^2\left(\Phi^{(t)}\right)
\eta^{-2} tr \left\{\left(
                                  \begin{array}{cc}
                                   - \eta_x\psi_y &  \frac{\eta_x\psi_x-\eta_y\psi_y}{2} \\
                                     \frac{\eta_x\psi_x-\eta_y\psi_y}{2} & \eta_y\psi_x \\
                                  \end{array}
                                \right) \left(
                             \begin{array}{cc}
                               |\Phi^{(t)}_x|^2 & \Phi^{(t)}_x\cdot\Phi^{(t)}_y \\
                               \Phi^{(t)}_x\cdot\Phi^{(t)}_y & |\Phi^{(t)}_y|^2  \\
                             \end{array}
                           \right)\right\}= \\
                           2\mu^2\left(\Phi^{(t)}\right)\eta^{-2}\left(
                          - \eta_x\psi_y  |\Phi^{(t)}_x|^2
                           +(\eta_x\psi_x-\eta_y\psi_y)\Phi^{(t)}_x\cdot\Phi^{(t)}_y
                           +\eta_y\psi_x|\Phi^{(t)}_y|^2 \right) \ .
                            \end{multline}
                            while
                           \be\label{sof}\mu^2(\Phi^{(t)})\eta^{-1}( \vec{v}\cdot
\nabla\eta)tr \left(|
D\Phi^{(t)}|^2\right)=\mu^2(\Phi^{(t)})\eta^{-2}\left(
\psi_y\eta_x-\psi_x\eta_y\right)\left( |\Phi^{(t)}_x|^2 +
|\Phi^{(t)}_y|^2\right) \ . \ee using (\ref{sqsi},\ref{crosspsi})
and (\ref{expand2**}-\ref{sof})and differentiation at $\tau=0$ yield
\begin{multline}\label{alfn}\frac{d}{dt} \hat{H}_\Delta\left(\Phi^{(t)}\right)=
\int_\Delta \mu^2\left(\Phi^{(t)}\right)\eta^{-1}
\left[\left(\Phi_x^{(t)}\cdot\Phi_y^{(t)}\right)\Box\psi +
\left(\left|\Phi_y^{(t)}\right|^2-\left|\Phi_x^{(t)}\right|^2\right)\psi_{xy}\right]dxdy
\\
+\int_\Delta \mu^2\left(\Phi^{(t)}\right)\eta^{-2}\left(
2(\eta_x\psi_x-\eta_y\psi_y)\Phi^{(t)}_x\cdot\Phi^{(t)}_y+
 (\eta_y\psi_x+\eta_x\psi_y)\left( |\Phi^{(t)}_y|^2
 -|\Phi^{(t)}_x|^2\right)\right)dxdy  \ .
\end{multline}
In complex notations (\ref{compnot}), (\ref{alfn}) takes the form
\begin{multline} \frac{d}{dt}\hat{H}_\Delta\left(\Phi^{(t)}\right)=  4
\mathfrak{Im}\left\{\int_\Delta
\mu^2\left(\Phi^{(t)}\right)\eta^{-1} f_{\Phi}\partial_{\oz}^2\psi
dxdy + \int_\Delta\mu^2\left(\Phi^{(t)}\right) (\partial_{\oz}
\eta^{-1})(\partial_{\oz} \psi)f_{\Phi}dxdy\right\}  \ \\
= 4\mathfrak{Im}\int_\Delta
\mu^2\left(\Phi^{(t)}\right)f_{\Phi}\partial_{\oz}\left( \eta^{-1}
\partial_{\oz}\psi\right)dxdy  \ .
\end{multline}
\end{proof}
\begin{proof} (of Proposition~\ref{nextDelta}): Integration by parts of (\ref{finalR}) yields:
\begin{multline}\label{finalRRR} \frac{d}{dt}\hat{H}_\Delta(\Phi^{(t)}) =
4\mathfrak{Im}\left(\int_\Omega\int \psi^{(t)}\partial_{\oz}\left\{
\partial_{\oz}\left( \mu^2(\Phi^{(t)}) f_{\Phi^{(t)}} \right)\eta^{-1}\right\}dxdy
\right) \\ -
2\mathfrak{Re}\left(\int_{\partial\Omega}\eta^{-1}\mu^2(\Phi^{(t)})f_{\Phi^{(t)}}
\partial_{\oz}\psi^{(t)}   dz\right) \ .
\end{multline}
The rest of the proof is, essentially, the same as that of
Proposition~\ref{nextDelta}.
\end{proof}
\begin{proof} (of Theorem~\ref{flowTh1} and Theorem~\ref{flowTh2}):
The proof of both Theorems follows from (\ref{finalRR}) and
(\ref{finalRRR}), respectively, where (\ref{finalbc}) is applied in
both cases.
\end{proof}

 \begin{center}{\bf References}\end{center}
 \begin{description}
 \item{[B]} \ Y. Brenier: {\it Polar factorization and monotone
rearrangement of vector valued  functions}, Comm. Pure Appl. Math,
{\bf 44}, (1991), 375-417.
 \item{[B2K]}
 A.M. Bronstein, M.M Bronstein, A.M Bruckstein and R. Kimmel, {\it
 Paretian similarity for partial compariaon of non-rigid objects},
 in  Scale space and variational methods in computer vision, F.
 Sgallari, A. Murli and N. Paragios (eds)., Springer 2007, 264-275.
\item{[C]} L Cafarelli: {\it  The Monge: Ampre equation and Optimal Transportation, an
Elementary Review},   2003 - Springer-Verlag  Lec. Notes, (2003)
 \item{[EK]} A. Elad and R. Kimmel,  {\it Bending invariants
 representations for surfaces}, Proc. CVPR, 2001, 168-174
\item{[FK]} H.M. Farkas and I Kra, {\it Riemann Surfaces},
Springer-Verlag, NY (1980)
 \item{[H]} F. H\'{e}lein, {\it Harmonic Maps, Conservation Laws and
 Moving Frames}, Cambridge Tracts in Mathematics {\bf 150}, Cambridge Univ.
 Press, 2002
 \item{[MS]} F.  M\'{e}moli and G. Sapiro, {\it A theoretical and
 computational framework for isometry invariant recognition of point
 cloud data}, Foundation of computational mathematics,  {\bf 5},
 2005, \#3, 313-347
 \item{[M]}  J.Moser, {\it On the volume elements on a manifold}
,"AMS Trans.120(2) (1965),pp.286-294
\end{description}
\end{document}